\documentclass[11pt]{amsart}
\usepackage{amssymb}

\newtheorem{theorem}{Theorem}[section]

\newtheorem{proposition}[theorem]{Proposition}

\newtheorem{lemma}[theorem]{Lemma}

\theoremstyle{remark}
\newtheorem{remark}[theorem]{Remark}

\newtheorem{example}[theorem]{Example}

\newcommand\be{\begin{equation}\label}
\newcommand\ee{\end{equation}}

\newcommand{\V}{\mathcal{V}}

\newcommand{\U}{\on{U}}

\newcommand{\R}{\mathbb{R}}
\newcommand{\C}{\mathbb{C}}

\newcommand{\Z}{\mathbb{Z}}

\newcommand\lie[1]{\mathfrak{#1}}

\newcommand{\h}{\lie{h}}

\newcommand{\g}{\lie{g}}

\newcommand{\n}{\lie{n}}
\newcommand{\z}{\lie{z}}

\newcommand{\on}{\operatorname}

\newcommand{\Ad}{ \on{Ad} }

\newcommand{\p}{\mf{p}}

\newcommand\beqn{\begin{equation}}      
\newcommand\eeqn{\end{equation}}

\newcommand{\mf}{\mathfrak}
\newcommand{\beq}{\begin{eqnarray*}}
\newcommand{\eeq}{\end{eqnarray*}}

\begin{document}

\title[Invariant *-products on coadjoint orbits]
{Invariant *-products on coadjoint orbits
and the Shapovalov pairing}

\author{A. Alekseev}
\address{University of Geneva, Section of Mathematics,
2-4 rue du Li\`evre, 1211 Gen\`eve 24, Switzerland}
\email{alekseev@math.unige.ch}

\author{A. Lachowska}
\address{Department of Mathematics, M.I.T., 77 Massachusetts ave., Cambridge 
MA 02139, USA}
\email{lachowska@math.mit.edu}

\date{\today}

\begin{abstract}
{We give an explicit formula for invariant $*$-products on a wide
class of coadjoint orbits. The answer is expressed in terms of the
Shapovalov pairing for generalized Verma modules.
}
\end{abstract}
\maketitle

\vskip 0.3cm

\section{Introduction}
The problem of constructing a $*$-product on a manifold
$M$ with given Poisson structure was formulated in \cite{5}.
In the case of symplectic manifolds the existence of
$*$-products was established in \cite{DL} and in a more geometric
fashion in \cite{F}. For general Poisson manifold the existence
of $*$-products was proved in \cite{K}. 
While general existence results are now available, giving
explicit formulas for $*$-products remains a difficult task.
The first formula of this type was given by Moyal in \cite{M}
in the case of a constant Poisson bi-vector on $\R^n$. 
Further examples that one can consider are linear Poisson
brackets on the dual $\g^*$ of the Lie algebra $\g$,
and nondegenerate Poisson brackets on coadjoint orbits
in $\g^*$. The $*$-products on $\g^*$ were constructed 
in {\em e.g.} \cite{G}. It is a natural idea  to construct 
$*$-products on coadjoit orbits by restriction from $\g^*$. 
However, in \cite{CGR} it was proved that for $\g$ semisimple 
a smooth $*$-product on $\g^*$ does not restrict to coadjoint orbits.

Examples of  $*$-products on  some simple coadjoint 
orbits ($C {\rm P}^n$ and symmetric spaces)
can be found {\em e.g.} in \cite{B}, \cite{T}.
We shall concentrate on constructing $*$-products
on $M$ {\em invariant} with respect to the transitive
$G$-action. In the case of $G={\rm GL}(n)$ this
problem has been addressed in \cite{DM1}.
Invariant $*$-products on minimal nilpotent
coadjoint orbits of simple Lie groups were 
constructed in \cite{ABC} and \cite{AB}. In these works
the locality axiom  (stating that the $*$-product should 
be defined by a bi-differential operator) is relaxed.
In the case of $\g$ semisimple constructions of 
$*$-products on coadjoint orbits of semisimple
elements were suggested in \cite{A} using the deformation
quantization with separation of variables
of \cite{Ka} and in \cite{DM2} using the
methods of the category theory. 

Our main result is an explicit formula for invariant
$*$-products on coadjoint orbits $G/H$ for which the
corresponding Lie algebras $\g$ and $\h$ (or their complexifications)
fit into a decomposition $\g=\n_+ \oplus \h \oplus \n_-$
(see Section \ref{sec:assume} for the precise statement
of assumptions). Examples include: the space $\R^{2n}$ with
constant nondegenerate Poisson bracket (a coadjoint
orbit in the dual to the Heisenberg algebra), coadjoint orbits
of semisimple elements in semisimple Lie algebras,
as well as infinite dimensional examples such as
coadjoint orbits in the dual to the Virasoro algebra. 
Our construction is motivated by the fusion techniques of \cite{EV}.

\vskip 0.2cm

{\bf Acknowledgements.} We received a lot of help from S. Parmentier who participated
in this project in its early stages. 
We are grateful to S. Gutt, P. Etingof, V. Rubtsov,
D. Sternheimer, X. Tang, B. Tsygan and P. Xu for the interest in our work and for useful
discussions. This work was supported in part by the Swiss National
Science foundation and by the Erwin Schr\"odinger Institute for
mathematical Physics. 

\section{Preliminaries}
In this Section we formulate the problem of
finding invariant $*$-products on homogeneous
spaces, and state the assumptions which allow
us to solve it in an explicit form. 

\subsection{Invariant *-products on homogeneous spaces}
Let $M$ be a manifold. Recall that a $*$-product on $M$
is defined by a formal power series in $\hbar$ with
coefficients complex bi-differential operators on $M$, 
$B= \sum_{n=0}^\infty \hbar^n B_n$ with $B_0=1$ 
and $B_n \in {\rm Diff}^2(M)$.
The $*$-product is then given by formula,
$f*g:= fg + \sum_{n=1}^\infty \hbar^n B_n(f,g)$
for $f,g\in C^\infty(M)$. The main condition imposed on 
$B$ is that the $*$-product be associative,
$f*(g*h)=(f*g)*h$ for all $f,g,h \in C^\infty(M)$.

Let $G$ be a connected Lie group and let $H\subset G$ be
a closed Lie subgroup of $G$. Denote the corresponding
Lie algebras by $\g$ and $\h$, respectively. The quotient
$M:= G/H$ carries a transitive action of $G$. There is an 
induced $G$-action on functions on $M$, and on differential and
poly-differential operators. A $*$-product on $M$ is called
invariant if it is defined by an invariant formal bi-differential
operator $B$. That is, all bi-differential operators $B_n$ 
have to be $G$-invariant.

Recall that the space of invariant differential operators
on $M=G/H$ can be expressed as follows, 
${\rm Diff}_G(M)=(U\g/U\g \cdot \h)^H$.
Here $U\g$ is the universal enveloping algebra of $\g$, 
$U\g \cdot \h$ is the left ideal generated by $\h \subset \g \subset U\g$.
The algebra $U\g$ carries a natural adjoint action of $H \subset G$,
and since $\Ad_H(\h) = \h$ this action factors to the
quotient $U\g/U\g \cdot \h$. The $H$-invariant part
$(U\g/U\g \cdot \h)^H$ has an algebra structure
induced by the one of $U\g$. In a similar fashion,
the space of invariant $N$-differential operators on $G/H$
is given by
$$
{\rm Diff}^N_G(M)=( (U\g/ U\g \cdot \h)^{\otimes N})^H,
$$
where the invariant part is taken with respect to the
diagonal $H$-action on $(U\g/ U\g \cdot \h)^{\otimes N}$.
Formal bi-differential operators
which define invariant $*$-products on $M$ take
values in the space
$$
B \in ((U\g/ U\g \cdot \h)^{\otimes 2})^H[[\hbar]] .
$$

\begin{remark}
As opposed to ${\rm Diff}_G(M)$, the space
${\rm Diff}^N_G(M)$ for $N > 1$ 
has no natural algebra structure.
\end{remark}

Let $\Delta: U\g \to U\g \otimes U\g$ be the standard coproduct
of $U\g$. That is, $\Delta$ is an algebra homomorphism such that 
$\Delta(x)=1\otimes x + x\otimes 1$ for all $x\in \g$.
Let $B,C \in {\rm Diff}^2_G(M)=( (U\g/ U\g \cdot \h)^{\otimes 2})^H$.
Then expressions $((\Delta \otimes 1)B)(C\otimes 1)$ and
$((1\otimes \Delta)B)(1\otimes C)$ define unique elements 
of ${\rm Diff}^3_G(M)=( (U\g/ U\g \cdot \h)^{\otimes 3})^H$.
In more detail, let $\hat{B}, \hat{C} \in U\g \otimes U\g$
be representatives of the classes 
$B,C\in (U\g/U\g \cdot \h)^{\otimes 2}$. The 
classes of $((\Delta \otimes 1)\hat{B})(\hat{C}\otimes 1)$ and
$((1\otimes \Delta)\hat{B})(1\otimes \hat{C})$ in 
$(U\g/U\g \cdot \h)^{\otimes 3}$ are $H$-invariant and
independent of the choice.

The associativity of an invariant $*$-product defined by
a formal bi-differential operator $B$ is expressed by the 
following equality of invariant formal 3-differential
operators,
\begin{equation} \label{eq:assoc}
\left( (\Delta \otimes 1) B \right) (B\otimes 1) =
\left( (1\otimes \Delta) B\right) (1\otimes B) 
\end{equation}
in $((U\g/U\g \cdot \h)^{\otimes 3})^H[[\hbar]]$.

\begin{remark}
If $H$ is a connected Lie group one can replace the condition
of $\Ad_H$-invariance of $B$ by an algebraic condition of
invariance with respect to the adjoint action of $\h$.
It is also a natural context when the Lie algebra $\g$ is infinite
dimensional and the corresponding Lie group may not be available.
\end{remark}

\begin{remark}
In the case of $\h=0$ equation \eqref{eq:assoc} was considered by
Drinfeld in \cite{Dr}. In this reference a family of solutions
of \eqref{eq:assoc} was constructed in term of the Campbell-Hausdorff
series for the Lie algebra $\g$. 
\end{remark}

\subsection{Assumptions on $\g$ and $\h$. Examples} \label{sec:assume}
Here we list assumptions imposed on
the Lie algebra $\g$ and its Lie sublagebra $\h$
which allow to construct explicit solutions of
equation \eqref{eq:assoc}. 

\begin{itemize}

\item The Lie algebra $\g$ is $\Z$-graded,
$ \g = \oplus_{i \in \Z} \g_i$, such that each 
graded component has finite dimension.
The adjoint action of $H$ on $\g$ preserves this
grading. 
Let $\h = \g_0$ and denote
$\n_+ = \oplus_{i >0}\g_i$, $\n_- = \oplus_{i<0}\g_i$.

\item There exists a character 
$\chi: \g_0 \to \C$ such that the paring $\n_+ \times \n_- \to \C$ 
defined by  $u,v \mapsto \chi([u,v]_0) $ for $u \in \n_+$, $v \in \n_-$ 
is nondegenerate.  Here $x \mapsto x_0$ denotes the projection 
onto the zero graded component $\g_0$. Such a character is
called nonsingular.

\end{itemize}

\begin{remark}
If $H$ is connected, the requirement that $H$ preserves
the grading is automatically satisfied.
\end{remark}

\begin{remark}
The second assumption implies that $\g_0 \neq 0$.
\end{remark}

\begin{example} \label{ex:quadr} 
Assume that $\g$ is quadratic. That is, there is a 
nondegenerate ${\rm ad}(\g)$-invariant symmetric bilinear
form $Q: \g \times \g \rightarrow {\C}$ of degree zero.
In this case, $\g_0$ is a finite dimensional quadratic 
Lie algebra with invariant bilinear form the restriction of
$Q$ to $\g_0$.
Let $\mathfrak{z}(\g_0)$ be the center of $\g_0$.
Every element $z \in \z(\g_0)$ defines
a character $\chi_z$ of $\g_0$,
$\chi_z(x) = Q(z, x)$.  An element 
$z \in \z(\g_0)$ is called nonsingular
if the map  $(u,v) \mapsto \chi_z([u,v]_0) \in \C$ 
defines a nondegenerate paring between the graded components 
$\g_i$ and $\g_{-i}$ for all $i$.
If $\z(\g_0)$  is nonempty and contains a nonsingular element, 
then the set of such elements is Zariski open in $\z(\g_0)$. 
\end{example}

\begin{example} \label{ex:ss1}
Let $\g$ be a semisimple complex Lie algebra,
$\h$ be a Cartan subalgebra of $\g$,
$\Pi =\{ \alpha_i \}_{i=1}^{{\rm rank}(\g)}$ be (some choice of) the set
of simple roots. The principal grading on $\g$
is the unique grading such that $\mathfrak{h} = \g_0$
and all root vectors corresponding to simple
roots have degree 1, $e_{\alpha_i} \in \g_1$.
Any regular character $\chi$ of $\mathfrak{h}$ defines a 
nondegenerate paring $\n_+ \times \n_- \to \C$. 

This construction applies verbatim to any Kac-Moody algebra.
\end{example}

\begin{example} \label{ex:ss2}

Let $\g$ be a semisimple complex Lie algebra, and let
$z \in \g$ be a semisimple element. Let $\g_z$ be 
the centralizer of $z$ and let $\mathfrak{h}$ be
a Cartan subalgebra of $\g$ such that
$z \in \mathfrak{h} \subset \g_z$. Define the 
unique grading on $\g$ such that $\g_0 = \g_z$
and for all simple roots in
$S_z:=\{ \alpha_i, e_{\alpha_i} \notin \g_z\}$
one has $e_{\alpha_i} \in \g_1$.
The Lie algebra $\mathfrak{p}_+=\g_0 \oplus \n_+$ is a
parabolic subalgebra of $\g$ with the Levi subalgebra $\g_z$,
corresponding to the subset $S_z\subset \Pi$. 
The center $\mathfrak{z}(\g_0) \subset
\mathfrak{h}$ consists of the elements orthogonal to $S_z$.
Denote by $R$ the set of roots of $\g$ and by $R_z \subset R$
the set of roots of $\g_z$. Then the nonsingular elements 
of $\mathfrak{z}(\g_0)$ that define nondegenerate bilinear 
paring between $\n_+$ and $\n_-$ are those with nonvanishing 
scalar product with all elements of $R\setminus R_z$.
They form a Zariski open subset of $\mathfrak{z}(\g_0)$.
In particular $z$ is a central nonsingular element in 
$\mathfrak{z}(\g_0)$.

If $z$ is a regular semisimple element of $\g$ we return
to Example \ref{ex:ss1}. 
\end{example}

\begin{example} \label{ex:Heis}
Let $\g = H_n$ be the Heisenberg
Lie algebra generated by $c$ and by $p_i, q_i, i=1 \dots n$ 
with the only nonvanishing Lie brackets $[p_i, q_j]=\delta_{ij} c$.
Define the grading by setting  ${\rm deg}(p_i)=1, {\rm deg}(q_i)=-1$ 
and ${\rm deg}(c)=0$. Then $\g_0 = \C c$ and any $\chi \in (\g_0)^*$ 
such that  $\chi(c) \neq 0$ is a nonsingular character. 
\end{example}

\begin{example} \label{ex:Vir}
Let $\g={\bf Vir} $ 
be the Virasoro algebra with $\g_n = \C L_n, n \in \Z, n \neq 0$ 
and  $\g_0 = \C L_0 \oplus \C \mf{c}$  with 
the Lie bracket $[L_n, L_m]= (n-m)L_{n+m} +\delta_{n+m, 0} 
\frac{n^3-n}{12}\mf{c}$.  Let $\chi \in \g_0^*$ be defined 
so that $\chi(L_0)= \Delta, \chi(\mf{c})=c$ for some $\Delta, c \in \C$. 
Then $\chi$ defines the paring  $(L_n, L_{-n})=2n\Delta + 
\frac{n^3-n}{12}c$ which is nondegenerate for all pairs 
$\Delta, c \in \C$ such that $\Delta + \frac{n^2-1}{24}c \neq 0$ for 
all $n \in \Z_{>0}$. 
\end{example}    

\subsection{Relation to coadjoint orbits}
Let $\g^*:=\oplus_{i\in \Z} \g_i^*$ be the 
(graded) dual to the Lie algebra $\g$. It carries
a natural coadjoint action of $\g$. An
element $\chi \in \g_0^* \subset \g^*$
defines a point in $\g^*$. Our second
assumption is equivalent to saying that
$\g_0$ is the coadjoint stabilizer of $\chi$.

Assume that $\g$ is finite dimensional. Then 
there is a connected simply connected Lie group 
$G$ with Lie algebra $\g$ and with a natural
coadjoint action of $G$ on $\g^*$. Denote
the stabilizer of $\chi$ by $G_0 \subset G$. It is a closed
subgroup of $G$ with Lie algebra $\g_0$. 
A coadjoint orbit of $\chi$ under the $G$-action
is a homogeneous space $M=G/G_0$. 
Our construction of invariant  $*$-products
applies if the adjoint action of $G_0$ preserves
the grading. If $G_0$ is connected, this condition
is satisfied automatically. In particular, this 
is always the case if $\g$ is a semisimple Lie algebra.

If $G_0$ is disconnected let $H$ be the connected component
of the unit element. Then  our method applies to $G/H$
which covers the coadjoint orbit $G/G_0$.

Let $G$ be a real Lie group and $\g$ be the corresponding
real Lie algebra. Since the $*$-products are usually defined
on the space of complex valued functions, it is sufficient that
assumptions of the previous Section be satisfied for the
complexified Lie algebra $\g^\C$. For this reason in the rest 
of the paper we assume that $\g$ is a complex Lie algebra.

Recall that the coadjoint orbits carry an invariant symplectic form
which is constructed as follows. Identify $T_{eH}(G/H)\cong \n_- \oplus \n_+$,
and define $\omega(u,v):= - \chi([u, v])$. By assumptions, this
bilinear form establishes a duality between $\n_-$ and $\n_+$.
Let $u_i$ and $v_i$ be a pair of dual bases in $\n_-$ and $\n_+$,
respectively. Then, the inverse of $\omega$ is the invariant
Kirillov-Kostant-Souriau Poisson bi-vector on $G/H$ which 
is equal to $\sum_i u_i \wedge v_i$ at
$eH$. An additional constraint which is often imposed on the 
bi-differential operator $B$ is that the skew-symmetric part
of $B_1$ be equal a given Poisson bi-vector. In our case, this
condition reads
\begin{equation} \label{eq:KKS}
B_1 - B_1^t = \sum_i u_i \wedge v_i ,
\end{equation}
where $B_1^t$ is a bi-differential operator obtained by 
exchanging two copies of $U\g/ U\g \cdot \h$ in the tensor product.
Geometrically, $\n_-$ and $\n_+$ define two distributions
in $T(G/H)$. If we deal with a real Lie group $\g$ and the conditions
of the previous Section are satisfied for the real Lie algebra $\g$,
these distributions give rise to transverse Lagrangian polarizations
on the orbit. In general, we are getting transverse complex 
polarizations. 

\section{Generalized Verma modules and the Shapovalov pairing}
In this Section we recall the notion and basic 
properties of the Shapovalov pairing and of the 
associated canonical element in the tensor product
of two opposite generalized Verma modules. 
The details on the generalized Verma modules 
can be found in e.g. \cite{Di}.

\subsection{Generalized Verma modules} 
Let $\p_=\oplus_{i \geq 0} \g_i$ and $\p_-=\oplus_{i \leq 0} \g_i$.
A central character $\chi$ of $\g_0$ can be given a 
$\p_\pm$-module structure by letting $\n_\pm$ act on it by zero.  
Define the generalized Verma modules by 
$$ M^\pm = {\rm Ind}_{U\p_\pm}^{U\g}  \chi 
\cong U(\g) \otimes_{U\p_\pm}  \chi.$$ 
In Example \ref{ex:ss1} we recover the Verma module over $\g$ 
of highest weight $\chi$, in Example \ref{ex:ss2} 
- the scalar generalized Verma module induced from the parabolic 
$\p_+$.  
By general properties 
of the induction, $M^+$ is isomorphic to $U\n_-$ 
as a $U\n_-$-module generated 
by $v_\chi = 1\otimes \chi$. 
Similarly, $M^- \cong U\n_+$ as a $U\n_+$-module generated by 
$v_\chi$. Both Verma modules inherit natural gradings 
from $U\n_+, U\n_-$. 

Suppose that $M$ is a $\Z$-graded $U(\g)$-module, and $V$ is any 
$U(\g)$-module. Then we define the completed tensor product
$M \check{\otimes} V \equiv \oplus_{i \in \Z} M_i \otimes V$
as a $\Z$-graded $U(\g)$-module 
where elements of $U(\g)$ act by comultiplication.  
In particular, if both 
$M$ and $N$ are $\Z$-graded $U(\g)$-modules, then 
$M \check{\otimes} N$ has two $\Z$-gradings. Often it is convenient 
to preserve both gradings,   
$ M \check{\otimes} N \equiv \oplus_{i,j}M_i \otimes N_j $, 
where the elements may have infinite length.

Fix a nonsingular character $\chi$ of $\g_0$. 
For any $\lambda \in \C$ we consider 
a rescaled character $\chi_\lambda: = \lambda \chi$, and
define a paring between 
$U\n_-$ and $U\n_+$ which depends on the parameter $\lambda$.
Write $U(\g) \cong U(\g_0) \oplus (\n_- U(\g) \oplus U(\g) \n_+)$ 
and let $\phi : U(\g) \rightarrow U(\g_0)$ be the projection of 
an element in $U(\g)$ to the first summand along the second. 
For any $x \in U\n_-, y \in U\n_+$, we set 
$$ (x,y)_\lambda := \chi_\lambda(\phi(S(y)x)), $$
where $S: U\g \to U\g$ is the antipode of $U\g$. That is,
$S$ is the unique anti-automorphism of $U\g$ such that
$S(x)=-x$ for all $x\in \g$.

Let 
$M^+_{\lambda} \cong {\rm Ind}_{U\p_+}^{U\g}  \chi_\lambda$ and 
$M^-_{-\lambda} \cong 
{\rm Ind}_{U\p_-}^{U\g}\chi_{-\lambda}$ be the generalized
Verma modules. 
The pairing $(\cdot, \cdot)_\lambda$ gives rise to a
pairing between $M^+_{\lambda}$ and $M^-_{-\lambda}$. 
Namely, fix generating vectors 
$v_\lambda = 1 \otimes \chi_\lambda \in M^+_\lambda$ and 
$v_{-\lambda} = 1 \otimes \chi_{-\lambda} \in M^-_{-\lambda}$
and let
$$( x v_\lambda , y v_{-\lambda}) := (x,y)_\lambda .$$ 
This pairing is called the Shapovalov pairing between
the generalized Verma modules. It is $U\g$-invariant in
the following sense: $(a u, v)= (u, S(a) v)$ for 
$u \in M^+_\lambda$ and $v \in M^-_{-\lambda}$. 
The modules $M^{+}_{ \lambda}$ and $M^-_{-\lambda}$
are irreducible if and only if 
$(\cdot, \cdot)_{\lambda}: U\n_+ \times U\n_- \to \C$ 
is nondegenerate. Indeed, if $x \cdot v_{\lambda}$ lies in a proper 
submodule of $M^+_{\lambda}$ and $x \in U\n_-$ has a maximal degree 
in this submodule,  then $(x, y)_{\lambda}=0$ for all $y \in U\n_+$ 
by the $U(\g)$-invariance of the paring, and conversely.

\begin{proposition} \label{prop:nonsing}
Let $\chi$ be a nonsingular character of $\g_0$. 
Then the pairing $(\cdot, \cdot)_\lambda: U\n_- \otimes U\n_+ \to \C$
is nonsingular for almost all $\lambda \in \C$.
\end{proposition}

\begin{proof}
Let $\{u_i\}$ be a homogeneous basis in $\n_-$, $\chi$ a nonsingular 
character, and $\{v_i\}$ the dual 
basis in $\n_+$ with respect to the paring 
$(u, v) = - \chi([u, v]_0)$. 
Choose an order 
of the elements $u_i$ in each graded component of $\n_+$, and 
enumerate the set $\{u_i\}$ by increased grading. 
Then we have a PBW-type basis in $U({\mathfrak n}_{\pm})$. In particular, 
in each graded 
component $(U\n_-)_{-n}, n \in \Z_+$, there is a basis  
$\{x^{(n)}_k\}_{k=1}^N$ of monomials in $\{u_i\}$. It can be  
ordered so that the number of factors in $x^{(n)}_i$ is greater or equal 
to that in $x^{(n)}_j$ whenever $i \leq j$. For each element 
$x^{(n)}_k = u_{k_1}^{s_1}u_{k_2}^{s_2} \ldots u_{k_r}^{s_r}$ set 
$y^{(n)}_k = v_{k_1}^{s_1}v_{k_2}^{s_2} \ldots v_{k_r}^{s_r}$.
The elements $\{y^{(n)}_k\}$ form a
basis in $(U\n_+)_{n}$. 
The number $d_k = \sum_{i=1}^r s_i$ will be called the length of $x^{(n)}_k$. 
Let $\lambda \in \C$ and 
consider the $N\times N$ matrix $M^n(\lambda)= (x^{(n)}_k, y^{(n)}_l)_{\lambda}$. 
Its elements 
are polynomials in $\lambda$.  
The following statements are easy to check: 

\begin{enumerate} 
\item{The order of the polynomial $(x^{(n)}_k, y^{(n)}_l)_{\lambda}$ 
cannot exceed the length of the shortest of the two 
monomials $x^{(n)}_k, y^{(n)}_l$.}
\item{For two monomials $x^{(n)}_k, y^{(n)}_l$ of the same length $d_k$, 
$(x^{(n)}_k, y^{(n)}_l)_{\lambda}$ is a polynomial of order strictly less than 
$d_k$ unless $l=k$. If $l=k$ then  $x^{(n)}_k$ and $y^{(n)}_k$ have the same number 
of factors which are dual to each other
with respect to the $(\cdot, \cdot)$ paring. }
\item{ We have $$ (x_k, y_k)_{\lambda} = (\prod_{i=1}^r (s_i)!) 
\lambda^{d_k} + P_k(\lambda),$$ 
where $P_k$ is a polynomial of order 
less than $d_k$.}   
\end{enumerate}
     
The above implies that we can write 
\begin{equation} \label{eq:M}
M^n(\lambda) = D^n(\lambda) C (1+O(1/\lambda)), 
\end{equation}
where $D^n(\lambda)$ is a diagonal matrix with 
$[D^n(\lambda)]_{kk}=(\prod_{i=1}^r (s_i)!) 
\lambda^{d_k} $, the matrix $C$ is  constant 
lower triangular with 
units on the diagonal, and $O(1/\lambda)$ is a matrix whose entries  
are polynomials in $1/\lambda$ without a constant term. 
The determinant of $M^n(\lambda)$ is a polynomial in $\lambda$ of 
order $\sum_{k=1}^N d_k$ with a nonzero leading coefficient. 
Therefore, the matrix $M^n(\lambda)$ is invertible for all but a 
finite number of values of $\lambda \in \C$. 
The union of zeros of $M^n(\lambda)$ for all $n$ is a countable 
subset of $\C$. This completes the proof.
\end{proof} 

\subsection{Canonical element $F_\lambda$}
Let $\lambda \in \C$ such that the pairing $(\cdot, \cdot)_\lambda$
be nonsingular. Denote by $F_\lambda \in U\n_- \check{\otimes} U\n_+$
the canonical element corresponding to the pairing. 
It has the following form, $F_\lambda \in 1 + (U\n_-)_{<0} \check{\otimes}
(U\n_+)_{>0}$. Let $M^+_\lambda \check{\otimes} M^-_{-\lambda}$
be a completed tensor product of two irreducible generalized
Verma modules with generating vectors $v_\lambda$ and $v_{-\lambda}$.
Then $F_\lambda(v_\lambda \otimes v_{-\lambda})$ is the canonical
element with respect to the Shapovalov pairing. In particular, it is
$U\g$-invariant. Choose a generating vector $v \in V_0$ in the trivial
$U\g$ module $V_0 \cong \C$. Then $v \mapsto F_\lambda(v_\lambda \otimes v_{-\lambda})$
defines a $U\g$-homomorphism $V_0 \to M^+_\lambda \check{\otimes} M^-_{-\lambda}$.

\begin{proposition} \label{prop:res}
The element $F_\lambda$ is a meromorphic function of $\lambda$
which is holomorphic at $\lambda =\infty$. The residue of
$F_\lambda$ at $\lambda=0$ is given by formula,
\begin{equation} \label{eq:res}
{\rm Res}_{\lambda=0} \, F_\lambda = \sum_{i=1}^n \, u_i \otimes v_i 
\in \n_- \check{\otimes} \n_+ \subset U\n_- \check{\otimes} U\n_+.
\end{equation}
\end{proposition}

\begin{proof}
Using bases $\{ x_k^{(n)} \}$ and $\{ y_k^{(n)} \}$ of Proposition
\ref{prop:nonsing} 
one can write the element $F_\lambda$ in the following form,
$$
F_\lambda = \sum_{n=0}^\infty \sum_{kl} (M^n(\lambda))^{-1}_{lk} \, 
x_k^{(n)} \otimes y_l^{(n)} .
$$
Decomposition \eqref{eq:M} implies that matrix elements  
$(M^n(\lambda))^{-1}_{lk}$ are rational functions of $\lambda$
holomorphic at $\lambda=\infty$. Hence, each bi-graded component
of $F_\lambda$ is a rational function of $\lambda$ holomorphic
at $\lambda=\infty$, which amounts to saying that $F_\lambda$
is meromorphic and regular at $\lambda=\infty$.

To compute the residue of $F_\lambda$ note that
for the matrix $C$ of \eqref{eq:M} one has
$C_{kl}=\delta_{kl}$  for $k,l$ with $d_k=d_l=1$.
This implies the same property for the inverse
matrix, $C^{-1}_{kl}=\delta_{kl}$ for $k,l$ with $d_k=d_l=1$.
We write,
$$
{\rm Res}_{\lambda=0}\, (M^n(\lambda))^{-1}_{lk}=
{\rm Res}_{\lambda=0}\, \sum_{m}(1+ O(\lambda^{-1}))_{lm}
C^{-1}_{mk} \lambda^{-d_k} (\prod_i s_i!)^{-1} =
{\rm Res}_{\lambda=0}\, C^{-1}_{lk} \lambda^{-d_k} (\prod_i s_i!)^{-1}
$$
This formula shows that the residue vanishes unless $d_k=1$.
Then, since $C^{-1}$ is lower triangular with $d_l\leq d_k$
for $l\geq k$ one has $d_l=1$ as well. Hence, the only
matrix elements which have nonvanishing residues
are diagonal entries $k=l$ with $d_k=1$, in which case
the residue is equal to one. Since the elements with
$d=1$ are basis elements in $\n_-$ and $\n_+$ we obtain
\eqref{eq:res} as required. 
\end{proof}

\begin{remark}  \label{rem:natural}
Formula \eqref{eq:M} implies that the coefficient
$(M^n(\lambda))_{lk}^{-1}$
in front of the element $x_k^{(n)} \otimes y_l^{(n)}$ in $F_\lambda$
is a rational function in $\lambda$ of order less or equal to  $- d_k$.
By exchanging the roles of $k$ and $l$ we see that 
in fact the order is $\leq - {\rm max}(d_k, d_l)$.
\end{remark}

\begin{remark}
In case when $\g$ is finite dimensional, let $d_{max} \in \Z_+$ be  
such that $(\n_+)_l =0$ for all $l >d_{max}$.
Then we have the following 
estimate for the length $d_k$ of monomials in 
the basis of $(U\n_+)_{n}$: $d_k \geq [\frac{n}{d_{max}}]$. 
Therefore in the Taylor series in $1/\lambda$ for the coefficients of
the matrix $(M^n(\lambda))^{-1}$ the leading term is of order greater or 
equal to $[\frac{n}{d_{max}}]$. The coefficients of a given order $m$ in $1/\lambda$ 
can appear only in matrices $(M^n(\lambda))^{-1}$ with $n \leq m d_{max}$. 
Since each graded component of $U\n_+$ is finite dimensional, $F_\lambda$ has  
only finitely many bi-graded components
with coefficients of a given order in $1/\lambda$.
\end{remark}

\section{Solutions of the associativity equation}

In this Section we relate the Shapovalov pairing 
with the associativity equation (\ref{eq:assoc}) for
the pair $\g, \h: = \g_0$.
In the proof of the associativity equation we shall use
$U\g$-homomorphisms from generalized Verma modules
to various completed tensor products. Our argument uses the method  
introduced in \cite{EV}; detailed exposition can be found in \cite{ES}.

\subsection{$U\g$-homomorphisms of generalized Verma modules}
One useful property of $U\g$-homomorphisms is given in the following
Proposition.

\begin{proposition} \label{prop:uniq}
Let $M^+_\lambda$ be an irreducible generalized Verma module 
with generating vector $v_\lambda$. Suppose that for a
$U\g$-module $\V$, there exists a $U\g$-invariant element 
$z \in M^+_\lambda \check{\otimes} \V$ such that $z \in v_\lambda 
\otimes w + (Un_-)_{< 0} \cdot v_\lambda \check{\otimes} \V$ 
for a certain  $w \in \V$. Then such  element is unique.  
\end{proposition}

\begin{proof}
Assuming that an invariant element exists, we will construct  it 
iductively starting from the summand $v_\lambda \otimes w$. 
Let $\{y_k\}$ and $\{x_k \}$ be the bases in $U\n_+$ and 
$U\n_-$ which consist of the elements of the bases  
in each graded component of $U\n_+$ and $U\n_-$
constructed in the proof of Proposition \ref{prop:nonsing}, 
and  ordered by the increased grading.  
Then we can write 
$z = v_\lambda \otimes w + \sum_{k>0}x_k v_\lambda \otimes w_k$ 
for some elements $w_k \in \V$.   
By assumptions  we have $\Delta(y_k) z=0$ for all $k>0$. 
Let $y_1$ be a basis element of degree $1$. Then 
$\Delta(y_1) z = v_\lambda \otimes y_1 w + \sum_k y_1 x_k v_\lambda 
\otimes w_k + \sum_k x_k v_\lambda \otimes y_1 w_k. $ 
By construction of the basis, there is only one $x_s$ (of degree 1) 
such that $y_1 x_s v_\lambda = v_\lambda$,  
and therefore one determines uniquely the element $w_s \in \V$. 
Then proceed by induction on 
the degree of $y_k$. 
\end{proof}

\begin{remark} 
In case when $\V=V$ is a finite dimensional module, $\check{\otimes}$ 
is the usual tensor product, the statement 
follows from the Frobenius reciprocity of the induction: 
$$ {\rm Hom}_{U\g}(V_0, M^+_\lambda \otimes V) = 
{\rm Hom}_{U\g}(M^-_{-\lambda},  V)= 
{\rm Hom}_{U\p_-}(\chi_{-\lambda}, V). $$
Therefore, the space of $U\g$-homomorphisms 
$ V_0 \to M^+_\lambda \otimes V$ 
coincides with the space of  
$(U\n_-)$-invariant vectors in $V$ with the action of $U\g_0$ given 
by the character $\chi_{-\lambda}$.
\end{remark}

\begin{proposition} \label{prop:A}
Let $A: \g \to \g$ be a Lie algebra automorphism of $\g$ preserving
the grading. Then the element $F_\lambda$ is invariant with respect to
the natural action of $A$ on $U\n_- \check{\otimes} U\n_+$.
\end{proposition}

\begin{proof}
Since $A$ preserves the grading, the element
$A(F_\lambda)(v_\lambda \otimes v_{-\lambda}) \in
M^+_\lambda \check{\otimes} M^-_{-\lambda}$ 
is of the form, $v_\lambda \otimes v_{-\lambda}
+ (U\n_-)_{<0} \cdot v_\lambda \otimes (U\n_+)_{>0} \cdot v_{-\lambda}$.
The map $A$ being a Lie algebra automorphism, 
$A(F_\lambda)(v_\lambda \otimes v_{-\lambda})$ is $U\g$-invariant.
Hence, by Proposition \ref{prop:uniq} one has
$A(F_\lambda)(v_\lambda \otimes v_{-\lambda})=
F_\lambda(v_\lambda \otimes v_{-\lambda})$ which implies 
$A(F_\lambda)=F_\lambda$.
\end{proof}

\begin{remark}
Assume that  $\g$ integrates to a Lie group $G$, and $H \subset G$
is a subgroup with Lie algebra $\g_0$ such that $\Ad_H$
preserves the grading. Then, the element $F_\lambda$
is $\Ad_H$-invariant.
\end{remark}

\begin{proposition} \label{thm:assoc}
Let $\lambda \in {\mathbb C}$ be such that the module 
$M^+_{\lambda}$ is irreducible and denote $p: U\g \to
U\g/U\g \cdot \g_0$ the natural projection. Then
\begin{equation} \label{eq:1p1}
(1\otimes p \otimes 1)\, 
\big[(\Delta \otimes 1)F_\lambda(F_\lambda \otimes 1) \big] =
(1\otimes p \otimes 1)\,
\big[ (1 \otimes \Delta) F_\lambda (1\otimes F_\lambda) \big]
\end{equation}
in $U\n_- \check{\otimes} U\g/U\g \cdot \g_0 \check{\otimes} U\n_+$. 
\end{proposition}

\begin{remark}
In case when $\g$ is semisimple and $\g_0$ 
its Cartan subalgebra, the equation (\ref{eq:1p1}) is a projection 
from $U(\g)^{\check{\otimes}3}$ to 
$U\n_- \check{\otimes} U\g/U\g \cdot \g_0 \check{\otimes} U\n_+$
of the dynamical cocycle condition which was considered in 
\cite{EV}, \cite{EE} and \cite{ESS}. The element satisfying such a 
condition is called a dynamical twist and is used furter to construct 
solutions of the dynamical Yang-Baxter equation. 

\end{remark}

For $V_0$ the trivial representation of $U\g$ define 
$\V_0 = {\rm Ind}_{U(\g_0)}^{U(\g)}V_0$. Choose a nonzero
vector $v \in V_0$ and denote $\tilde{v}:= 1\otimes v$
the generating vector in $\V_0$.

\begin{lemma} \label{lem:F+-}
Let $\lambda \in {\mathbb C}$ be such that the module 
$M^+_{\lambda }$ is irreducible.
\begin{enumerate}
\item{ There exists a 
unique $U(\g)$-homomorphism 
$$ {\mathcal F}_\lambda^+ : 
M^+_{\lambda } \to M^+_{\lambda } 
\check{\otimes} \V_0 $$ 
such that 
$${\mathcal F}_\lambda ^+ (v_{\lambda }) \in  
v_{\lambda } \otimes \tilde{v} +  (U\n_-)_{<0} \cdot v_{\lambda }\otimes 
U\g \cdot \tilde{v}.$$
Explicitly,
$$ {\mathcal F}_\lambda ^+ (v_{\lambda }) = 
F_\lambda(v_{\lambda }\otimes \tilde{v}).$$}
\item{Similarly, there exists a 
unique $U\g$-homomorphism 
$$ {\mathcal F}_{-\lambda}^- : M^-_{-\lambda } \to \V_0 
\check{\otimes}M^-_{-\lambda } $$ 
such that for a fixed vector $v \in V_0$
$$ {\mathcal F}_\lambda ^- (v_{-\lambda } ) \in  
\tilde{v} \otimes v_{-\lambda } + 
U\g \cdot \tilde{v} \otimes (U\n_+)_{>0} \cdot v_{-\lambda }.$$
Explicitly,
$$  {\mathcal F}_{-\lambda}^- (v_{-\lambda }) = F_\lambda 
(\tilde{v} \otimes v_{-\lambda }).$$}
\end{enumerate}
\end{lemma}

\begin{proof}
Consider the first statement. 
The element $F_\lambda(v_{\lambda} \otimes \tilde{v})$ is $U\n_+$-invariant 
because of the $U\g$-invariance of the Shapovalov's paring and 
the $U\n_+$-invariance of $v_\lambda$. The subalgebra $U\g_0$ acts 
on $F_\lambda(v_{\lambda} \otimes \tilde{v})$ 
by the character $\chi_\lambda$. Therefore, by the universal 
property of the generalized Verma module, there is a homomorphism 
${\mathcal F}^+_\lambda: M^+_{\lambda } \to M^+_{\lambda } 
\check{\otimes} \V_0 $ mapping the generating vector $v_\lambda$
to $F_\lambda(v_{\lambda} \otimes \tilde{v})$.
The proof of the uniqueness for a given choice of $v \in V_0$ coincides 
verbatim with the proof of Proposition \ref{prop:uniq}. 
The second statement is proved similarly.   
\end{proof}

\begin{lemma} \label{lem:two}
Let $\lambda \in {\mathbb C}$ be such that $M^+_{\lambda}$
is irreducible. 
The two morphisms of $U\g$-modules 
$$ V_0 \rightarrow
M^+_{\lambda} \check{\otimes} M^-_{-\lambda} 
\stackrel{{\mathcal F}^+_\lambda \check{\otimes} 
{\rm id}}{\longrightarrow} 
M^+_{\lambda} \check{\otimes} \V_0 \otimes M^-_{-\lambda}$$ 
and 
$$ V_0 \rightarrow
M^+_{\lambda} \check{\otimes} M^-_{-\lambda} 
\stackrel{{\rm id} \otimes {\mathcal F}^-_\lambda}{\longrightarrow} 
M^+_{\lambda} \check{\otimes} \V_0 \check{\otimes} M^-_{-\lambda}$$ 
coincide. 
\end{lemma}

\begin{proof}
Both homomorphisms map $v$ to 
$v_\lambda \otimes \tilde{v} \otimes v_{-\lambda} + 
(U\n_-)_{<0}\cdot v_{\lambda} 
\check{\otimes} U(\g) \cdot \tilde{v} \check{\otimes} 
(U\n_+)_{>0} \cdot v_{-\lambda}$. 
Such a homomorphism is unique by Lemma \ref{lem:F+-}.     
\end{proof}

\begin{proof}[Proof of Proposition \ref{thm:assoc}]
Written in terms of 
$F_\lambda \in \sum_i (U\n_+)_{-i} \otimes 
(U\n_-)_i$, the equality of the two homomorphisms
in Lemma \ref{lem:two} gives the statement of the theorem.  
\end{proof}

\subsection{Invariant $*$-products}
Denote by $\pi: U\g^{\otimes 2} \to (U\g/U\g \cdot \g_0)^{\otimes 2}$
the natural projection. The following theorem is the main result
of this paper.

\begin{theorem}
Let $\chi$ be a nonsingular character of $\g_0$ and
$F_\lambda$ be the corresponding canonical element
in $U\n_- \check{\otimes} U\n_+$. Then,
the element $B:=\pi(F_{\hbar^{-1}})$ takes values in
$((U\g/U\g \cdot \g_0)^{\otimes 2})^{\g_0}[[\hbar]]$ and
satisfies the associativity equation \eqref{eq:assoc}. 
\end{theorem}

\begin{proof}
Recall that $F_\lambda$ is a meromorphic function of $\lambda$
holomorphic at $\lambda =\infty$. Hence, $F_{\mu^{-1}}$
is holomorphic at zero and defines a formal power series
$F_{\hbar^{-1}} \in (U\n_- \check{\otimes} U\n_+)[[\hbar]]$
with projection $B:=\pi(F_{\hbar^{-1}}) \in 
(U\g/U\g \cdot \g_0)^{\otimes 2}[[\hbar]]$. The $U\g$-invariance
of the Shapovalov pairing implies the $\g_0$-invariance
of $\pi(F_\lambda)$ and as a consequence the $\g_0$-invariance
of $B$. 

Projecting equation \eqref{eq:1p1} to $(U\g/U\g \cdot \g_0)^{\otimes 3}$
yields the associativity equation for $\pi(F_\lambda)$. Since
the latter is holomorphic at infinity, the associativity equation for
$B$ follows.
\end{proof}

\begin{remark}
As usual, assume that  $\g$ integrates to a Lie group $G$ and $\g_0$
integrates to a subgroup $H\subset G$ such that $\Ad_H$
preserves the grading. Then $\pi(F_\lambda)$ is an element in 
$((U\g/U\g \cdot \g_0)^{\otimes 2})^{H}[[\hbar]]$ which satisfies
the associativity equation. Hence, it defines an invariant
$*$-product on $G/H$.
\end{remark}

\begin{remark}
In a recent paper \cite{DM2}, a categorical approach to constructing 
dynamical twists is developed. It allows the authors 
to obtain a quantization of function algebras on semisimple 
coadjoint orbits for $\g$ reductive and $\g_0$ its Levi subalgebra.
The difference between the approach
of \cite{DM2} and the one chosen in the present paper is that
we do not use finite dimensional representations of $U\g$ and
the harmonic analysis on $G/H$.  
\end{remark}

\begin{proposition}
Let $\{ u_i\}$ and $\{v_i \}$ be dual bases in $\n_-$ and $\n_+$
with respect to the pairing $-\chi([\cdot, \cdot]_0)$. Then the element
$B$ has the form
$$
B=1 + \hbar \pi(\sum_i u_i \otimes v_i)  + O(\hbar^2).
$$
\end{proposition}

\begin{proof}
By Proposition \ref{prop:res} the first two terms in the
Taylor expansion of $F_\lambda$ at $\lambda=\infty$ are as
follows,
$$
F_\lambda = 1 + \lambda^{-1} \sum_i u_i \otimes v_i + O(\lambda^{-2}).
$$
Replacing $\lambda^{-1} \mapsto \hbar$ and applying projection $\pi$
yields the result.
\end{proof}

\begin{remark}
The skew-symmetric part of the first order term in $\hbar$ in the
Taylor expansion of $B$ defines an invariant bi-vector on
$M=G/H$. Its value at $eH \in M$ is given by $\sum_{i=1}^n u_i \wedge v_i$
which is exactly the Kirillov-Kostant-Souriau bi-vector 
on the coadjoint orbit. 
\end{remark}

\begin{remark} 
By Remark \ref{rem:natural} the elements $x_k \otimes y_l$ which
occur with a factor $\lambda^{-n}$ in the element $F_\lambda$
have the property $d_k, d_l \leq n$. Hence, the bi-differential
operators $B_n$ in the $*$-product $\pi(F_{\hbar^{-1}})$ have the order
$\leq n$ in each factor. Such $*$-products are called {\em natural}.
In \cite{GR} it is shown that a natural $*$-product induces a symplectic
connection on the underlying manifold $M$. In the case of coadjoint
orbits considered in this paper such a connection is defined
by a $\g_0$-invariant complement to to $\g_0$ in $\g$,
$\n_+ \oplus \n_- \subset \g$. 
\end{remark}

\begin{example}
Let $\g= H_n$ be the Heisenberg algebra defined in Example \ref{ex:Heis}, and 
let $\chi : c \mapsto w \in \C$ be a nonsingular character, $w \neq 0$.
Then the matrix of the Shapovalov paring in each graded component 
$(\U(n_-)_{-k}, \U(n_+)_k)$ is diagonal, with matrix elements
$$
(q_1^{k_1} \dots q_n^{k_n}, p_1^{k_1} \dots p_n^{k_n})_\lambda=
(-\lambda w)^k \prod_{i=1}^n k_i! ,
$$
where $k=k_1 + \dots + k_n$.
The corresponding inverse element $F_\lambda$ is given by
$$
F_\lambda = \sum_{k_1, \dots k_n=0}^\infty \,
\frac{(-1)^k}{(\lambda w)^k k_1! \dots k_n!} \,  
q_1^{k_1} \dots q_n^{k_n} \otimes p_1^{k_1} \dots p_n^{k_n} .
$$
Using the product structure of $U\n_- \otimes U\n_+$ one
can write the answer for $B$ in the following compact form,
$$
B=\exp\left(-\frac{\hbar}{w} \, \sum_{i=1}^n \, q_i \otimes p_i\right) .
$$
This gives a $*$-product on the coadjoint orbit $G/G_0 = \R^{2n}$ 
of the Heisenberg group acting on $\g^*$.
It is a `normal ordered' version of the Moyal product \cite{M}.
\end{example}

\begin{example}
Let $\g =sl(2, \C)$ with generators $e,f,h$ and the Lie brackets 
$[e,f]=h, [h,e]=2e, [h,f]=-2f$. Then $\g_0 = \C h$ and $\chi(h)=z \in \C$, 
$z \neq 0$ defines a nonsingular character.  The element $B$ associated
to $sl(2, \C)$ and $\chi$ is given by the formula 
\begin{equation} \label{eq:sl2}
B = \sum_{n=1}^\infty \hbar^n \left[
\frac{(-1)^n}{n!z(z-\hbar)\ldots (z-(n-1)\hbar)} \right]f^n \otimes e^n.
\end{equation}

\begin{remark}
By the theorem of Cahen, Gutt and Rawnsley \cite{CGR}, in case when $\g$ 
is semisimple, there are no $*$-products on $\g^*$ with the first order
term given by the KKS Poisson bi-vector 
which restrict to the coadjoint 
orbits. This is illustrated by the above example: starting from 
the second graded component, the expression \eqref{eq:sl2} for the $*$-product 
is singular at $z =0$. 
\end{remark}

\begin{remark}
Set $G= SU(2)$. Then the stabilizer of a nonsingular character $\chi \in \g^*$ is 
isomorphic to $G_0 =U(1)$. The homogeneous space $M = 
G/G_0 = S^2$ is the underlying manifold for the $*$-product constructed 
above.  
In this case explicit formula for the $*$-product similar
to \eqref{eq:sl2} appeared before 
in physics literature in \cite{P}, \cite{N} and \cite{HNT}.
\end{remark}

\end{example}

\begin{example}
Let ${\bf Vir}$ be the Virasoro algebra defined in Example \ref{ex:Vir} and 
$\chi$ a nonsingular character.  We will use the notation of  
\ref{ex:Vir} 
to give an explicit expression for the element $B$, corresponding to 
${\bf Vir}$ and  $\chi$, 
up to the second graded component. The character   
$\chi$ being nonsingular implies in particular that $\Delta \neq 0$ and 
$A := -32\Delta^3 -4 \Delta^2c \neq 0$. Denote $B := 20 \Delta^2 -2 \Delta c$, 
and set $D = A +\hbar B$. Then 
\begin{equation*} \begin{array}{l}  
 B = 1 - \hbar \left[\frac{1}{2\Delta}\right]L_{-1}\otimes L_1 + 
\hbar \left[\frac{8\Delta^2}{D}+\hbar\frac{4\Delta}{D}\right]
L_{-2}\otimes L_2 + 
 \hbar^2 \left[\frac{6\Delta}{D}\right]L_{-2}\otimes L_1^2 - \\
- \hbar^2 \left[\frac{6\Delta}{D}\right]L_{-1}^2 \otimes L_2 +
\hbar^2 \left[\frac{A}{8D}\right] L_{-1}^2 \otimes L_1^2 + 
(U\n_-)_{\leq -3} \otimes (U\n_+)_{\geq 3} .
\end{array} 
\end{equation*}

\end{example}

\end{document}